\documentclass{sigma}

\def\mfS{{\mathfrak S}}
\def\c{{\mathfrak c}}
\def\p{{\mathfrak p}}

\def\N{{\mathbb N}}

\def\X{{\mathbb X}}
\def\Y{{\mathbb Y}}

\def\x{\mathbf{x}}
\def\y{\mathbf{y}}
\def\z{\mathbf{z}}

\def\cE{{\mathcal E}}

\def\moins{\raise 1pt\hbox{{$\scriptstyle -$}}}
\def\plus{\raise 1pt\hbox{{$\scriptstyle +$}} }

\def\CARRE#1{\hbox{\vrule width \thickness
   \vbox to \carresize{\hrule height \thickness\vss
      \hbox to \carresize{\hss#1\hss}
   \vss\hrule height\thickness}
\unskip\vrule width \thickness} \kern-\thickness}
\def\vsquare#1{\vbox{\CARRE{$#1$}}\kern-\thickness}
\def\blk{\omit\hskip\carresize}

\def\young#1{%
  \newdimen\carresize \carresize=16pt%
  \newdimen\thickness \thickness=0.5pt%
  \vcenter{%
    \vbox{\smallskip\offinterlineskip%
      \halign{&\vsquare{##}\cr #1}}}}

\begin{document}
\allowdisplaybreaks

\renewcommand{\PaperNumber}{029}

\FirstPageHeading

\renewcommand{\thefootnote}{$\star$}

\ShortArticleName{The 6 Vertex Model and Schubert Polynomials}

\ArticleName{The 6 Vertex Model and Schubert
Polynomials\footnote{This paper is a contribution to the Vadim
Kuznetsov Memorial Issue `Integrable Systems and Related Topics'.
The full collection is available at
\href{http://www.emis.de/journals/SIGMA/kuznetsov.html}{http://www.emis.de/journals/SIGMA/kuznetsov.html}}}

\Author{Alain LASCOUX}

\AuthorNameForHeading{A. Lascoux}

\Address{Universit\'e de Marne-La-Vall\'ee, 77454,
Marne-La-Vall\'ee, France}
\Email{\href{mailto:Alain.Lascoux@univ-mlv.fr}{Alain.Lascoux@univ-mlv.fr}}
\URLaddress{\url{http://phalanstere.univ-mlv.fr/~al/}}

\ArticleDates{Received October 24, 2006; Published online February
23, 2007}

\Abstract{We enumerate staircases with f\/ixed left and right
columns. These objects cor\-respond to
 ice-conf\/igurations, or alternating sign matrices,
with f\/ixed top and bottom parts. The resulting partition
functions are equal, up to a normalization factor, to  some
Schubert polynomials.}


\Classification{05E15; 82B23}

\begin{flushright}\it
To the memory of Vadim Kuznetsov
\end{flushright}

\section{Introduction}

The $6$ vertex model is supposed to tell something about ice,
though not of the kind used in cold drinks. One has a planar array
f\/illed with  water molecules having $6$ possible types of
orientation. One chooses $6$ functions of the two coordinates.
Each molecule is weighted by one these functions, according to its
type. The weight of a given `square ice' is the product of the
weights of the molecules composing it.

In our case, the weights will be $x_i/y_j$, or $x_i/y_j -1$ for
two types of molecules, $i$, $j$ being the coordinates, the other
four types having weight $1$.

There are many other combinatorial objects equivalent to
square-ice conf\/igurations: alternating sign matrices (ASM),
totally symmetric self-complementary plane partitions(TSSCPP),
monotone triangles, staircases (subfamily of Young tableaux). We
choose this last object, because it is the most compact, and
because we are going to relate the combinatorics of weights to the
usual combinatorics of tableaux. We shall refer to Bressoud
\cite{Bressoud} for a description of the dif\/ferent approaches to
square-ice conf\/igurations.

There is a huge literature concerned with the enumeration of plane
partitions, and of  dif\/ferent related combinatorial objects,
starting from the work of MacMahon (the book of
Bressoud~\cite{Bressoud} provides many references). We are only
concerned here with the enumeration of ASM's, which started with
the work of Robbins and Rumsey \cite{MRR,RR}, with notable
contributions of Zeilberger~\cite{Zeilberger}
 and Kuperberg \cite{Kuperberg}.

One can put dif\/ferent weights on ice conf\/igurations or ASM.
The weight chosen by physicists  allows, by specialization, the
enumeration of ASM, after the evaluation of the Izergin--Korepin
determinant \cite{Izergin}. Kirillov and Smirnov \cite{Kirillov}
obtained determinants which correspond to more general partition
functions. Gaudin \cite[Appendix B, p.72]{Gaudin} had previously
given for the Bethe model a determinant similar to the
Izergin--Korepin determinant.

Okada \cite{Okada}, Hamel--King \cite{HK,HK2} obtained sums over
sets of ASM, with another weight involving  one coordinate only.
Their formulas generalize the interpretation of a Schur function
as a sum of Young tableaux of a given shape.

The weights that we have chosen involve the two coordinates, and
generalize the \emph{inversions} of a permutation; we  already
used them to obtain \emph{Grothendieck polynomials} from sets of
ASM~\cite{LaXGASM}. The sets that we take in this text are
simpler. They are sets of ASM having top and bottom part f\/ixed
(the case where one f\/ixes the top row an ASM already appears in
\cite{MRR,RR}). For these sets the `partition function' is
proportional to some Schubert polynomial, or to some determinant
generalizing Schur functions.

One remarkable feature of the functions that we obtain is a
symmetry property in the $x$-va\-riab\-les. This property is
dif\/ferent from the symmetry related to the Yang-Baxter equation
that the Izergin--Korepin determinant displays (see also
\cite{La99}).

The weight that we have chosen cannot be specialized to $1$ (and
thus, does not allow plain enumeration of ASM). On the other hand,
specializing all the $x_i$'s to $2$, and $y_i$'s
 to $1$ amounts to  weigh an ASM by $2^k$,
where $k$ is the number of $\moins 1$ entries. The `2'-enumeration
of all ASM, or to ASM having a top row f\/ixed,
 is due to \cite{MRR}, our study gives it
for the ASM having f\/ixed top and bottom parts of any size.

\section{Reminder}

\emph{Divided differences} $\partial_\sigma$ are operators on
polynomials in $x_1,x_2,\ldots$ indexed by permutations $\sigma$.
They are products of \emph{Newton divided differences}
$\partial_i:=\partial_{s_i}$ corresponding to the case of simple
transpositions $s_i$, and acting as
\begin{gather*}
\partial_i\colon f \to \frac{ f - f^{s_i}}{x_i-x_{i+1}}
\end{gather*}
denoting by $f^{s_i}$ the image of $f$ under the exchange of $x_i$
and $x_{i+1}$.

In the case where $\sigma=\omega := [n,\ldots, 1]$ is the maximal
permutation of the symmetric group~$\mfS_n$ on $n$ letters, then
$\partial_\omega$, apart from being expressible as a product of
$\partial_i$'s, can also be written
\begin{gather*}
f \to  \sum_{w\in \mfS_n} (-1)^{\ell(w))} f^w \prod_{1\leq i<j\leq
n} (x_i-x_j)^{-1},
\end{gather*}
where $\ell(w)$ is the \emph{length} of $w$.

The \emph{code} $\c(\sigma)$ of the permutation $\sigma$ is the
sequence of numbers $c_i:= \#(j\colon   j>i, \sigma_j<\sigma_i)$.

\emph{Schubert polynomials} are polynomials in two sets of
variables $\x$, $\y$ indexed by permutations or by their codes. We
use the symbol $\X$ in the f\/irst case, $\Y$ in  the second:
\begin{gather*}
\X_\sigma( \x, \y) = \Y_{\c(\sigma)} (\x,\y).
\end{gather*}

A \emph{dominant permutation} is a permutation with weakly
decreasing code. A \emph{Gra{\ss}mannian permutation} is a
permutation with only one descent, i.e.\ such that there exists
$r$:
\begin{gather*}
\sigma_1<\cdots < \sigma_r, \quad \sigma_{r+1}<\sigma_{r+2}, \quad
\dots
\end{gather*}
($r$ is the \emph{descent} of $\sigma$).

\emph{Dominant Schubert polynomials} $\Y_v(\x,\y)$,
 ( for $v_1\geq v_2 \geq v_3\geq \cdots$), are equal to
\begin{gather*}
\Y_v(\x,\y) := \prod_i \prod_{j=1}^{v_i}  (x_i-y_j).
\end{gather*}

General Schubert polynomials are by def\/inition all the non-zero
images of dominant Schubert polynomials under divided
dif\/ferences \cite{CBMS,MacSchub}.

The Gra{\ss}mannian  Schubert polynomials $\X_\sigma(\x,\y)$,
$\sigma$ having descent in $r$, are the images under
$\partial_{\omega_r}$ of the dominant polynomials $\Y_v$, $v\in
\N^r$. They have a determinantal expression. Let
$v=[\sigma_1-1,\ldots, \sigma_r-1]$. Then
\begin{gather*}
\X_\sigma(\x,\y) = S_{v_1-0,\ldots, v_r-r+1}(\x^r-\y^{v_1},\ldots,
\x^r-\y^{v_r}) = \det  \bigl|  S_{v_j+1-i}( \x^r - \y^{v_j})
\bigr|_{i,j=1,\ldots,r},
\end{gather*}
where $\x^r=[x_1,\ldots, x_r]$, $\y^k=[y_1,\ldots,y_k]$, and the
\emph{complete function}
 $S_m(\x^r -\y^k)$ is def\/ined as the coef\/f\/icient of $z^m$ in
\begin{gather*}
\prod_{1}^{k} (1-z y_i)   \prod_{1}^{r} (1-z x_i)^{-1}.
\end{gather*}

Multiplication of such a polynomial by $x_1\cdots x_r$ is easy,
starting from the appropriate dominant polynomial. It reduces to a
uniform increase of the index of the Schur function:
\begin{gather}
\label{ShiftGrassm} x_1\cdots x_r  \X_\sigma(\x,\y) =
 S_{v_1+1,\ldots, v_r-r+2}(\x^r-\y^{v_1},\ldots, \x^r-\y^{v_r}).
\end{gather}

The Gra{\ss}mannian Schubert polynomials specialized to
$\y=\mathbf{0}$ are exactly the  Schur functions. When specialized
to $\y=[0,1,2,\ldots]$, they coincide with the \emph{factorial
Schur functions}.

We shall need more general Schur-like functions.

Given $u,v\in \N^n$, and $2n+1$ alphabets $\z,\y_1,\ldots, \y_n,
\x_1,\ldots, \x_n$, then
\begin{gather*}
 S_{v/u}(\z-\y_1,\ldots, \z-\y_n // \x_1,\ldots, \x_n)
:= \det  \bigl| S_{v_j-u_i+j-i} (\z+\x_i - \y_j) \bigr|.
  \end{gather*}

Notice that if $\x_n, \x_{n-1},\ldots,\x_1$ are of cardinality
majorized respectively by $0,1,\ldots, n-1$, then  \cite[Lemma
1.4.1]{CBMS}
\begin{gather}
  S_{v/u}(\z-\y_1,\ldots, \z-\y_n // \x_1,\ldots, \x_n)
  = S_{v/u}(\z-\y_1,\ldots, \z-\y_n)\nonumber
\\ \qquad
{} :=  S_{v/u}(\z-\y_1,\ldots, \z-\y_n //
\mathbf{0},\ldots,\mathbf{0}). \label{SchurDrap}
\end{gather}

\section{Staircases}

A \emph{column} $u$ is a strictly decreasing sequence of integers,
its length is denoted $\ell(u)$. A \emph{staircase} is a sequence
of columns of lengths $k,k-1,\ldots,k-r$, such that, writing them
in the Cartesian plane as a Young tableau, aligning their bottoms,
then rows are weakly increasing, diagonals are weakly decreasing
(this last condition is added to the usual def\/inition of a Young
tableau, \emph{staircases are special Young tableaux}):
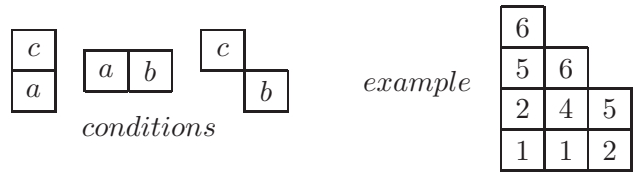
\begin{gather*}
\begin{array}{c}
  \young{c \cr a\cr} \quad \young{ a & b\cr} \quad \young{c\cr \blk & b\cr} \\
 conditions \end{array}
\,   \qquad example\quad
  \young{ 6\cr 5 &6\cr 2&4&5\cr 1&1 &2\cr}\,
\end{gather*}
with $c>a$, $a\leq b\leq c$.

Staircases are fundamental in the description of the
Ehresmann--Bruhat order of the symmetric group \cite{LS}.

Given two columns $u$, $v$ with $\ell(u)>\ell(v)$, let $\cE(u,v)$
be the set of all staircases with f\/irst column $u$, last column
$v$. We write   $\cE(n,v)$ when $u=[n,\ldots,1]$.

Given a staircase written in the plane as a tableau $t$, with
entries denoted $t[i,j]$, we give a~weight to each entry of $t$ as
follows:
\begin{itemize}\itemsep=0pt
\item  $t[i,j]=b$ has weight $x_i y^{-b} -1$ if $t[i-1,j] =b$;
\item  $t[i,j]=b$ has weight $x_i y^{-b}$ if  $t[i-1,j] < b<
t[i-1,j+1]$; \item in all other cases, $t[i,j]$ has weight $1$;
\item the weight of the staircase is the product of all these
elementary weights.
\end{itemize}

In other words, entries in the f\/irst column have weight $1$, and
there are three possible conf\/igurations for the other entries,
according to their left neighbors, with $a<b<c$:
\begin{gather*}
\begin{array}{lccc}
  &  \begin{matrix} c\\ b & b   \end{matrix}
  & \quad \begin{matrix} c\\ a & b   \end{matrix}
  & \qquad  \begin{matrix} b\\ a & b   \end{matrix}  \\
 weight   & \dfrac{x_i}{y_b}-1  &\ \  \dfrac{x_i}{y_b}    &\quad  1
\end{array}
\end{gather*}

Given two columns $u,v$, the sum of weights of the staircases in
$\cE(u,v)$ (resp.~$\cE(n,v)$) is denoted $F(u,v)$
(resp.~$F(n,v)$).  We call it the \emph{partition function} of the
set of staircases. When taking the alphabet $z_1,z_2,\ldots$
instead of $x_1,x_2,\ldots$, we write $F(u,v; z_1,z_2,\ldots)$.

Given an ASM, one builds a new matrix by replacing each row by the
sum of all rows above. The successive columns of the staircase
record the positions of the $1$'s in the successive rows of this
new matrix:
\begin{gather*}
\begin{bmatrix}
0   &   0   &   0   &   1   &   0   &   0   \\[.4ex]
0   &   0   &   1   &   -1   &   1   &   0   \\[.4ex]
0   &   1   &   -1   &   1   &   -1   &   1   \\[.4ex]
1   &   -1   &   1   &   -1   &   1   &   0   \\[.4ex]
0   &   1   &   -1   &   1   &   0   &   0   \\[.4ex]
0   &   0   &   1   &   0   &   0   &   0
\end{bmatrix} \quad
\begin{bmatrix}
0   &   0   &   0   &   1   &   0   &   0   \\[.4ex]
0   &   0   &   1   &   0   &   1   &   0   \\[.4ex]
0   &   1   &   0   &   1   &   0   &   1   \\[.4ex]
1   &   0   &   1   &   0   &   1   &   1   \\[.4ex]
1   &   1   &   0   &   1   &   1   &   1   \\[.4ex]
1   &   1   &   1   &   1   &   1   &   1
\end{bmatrix}   \quad
\young{
 6 \cr
 5 & 6 \cr
 4 & 5 & 6 \cr
 3 & 4 & 5 & 6 \cr
 2 & 2 & 3 & 4 & 5 \cr
 1 & 1 & 1 & 2 & 3 & 4 \cr
}\ .
\end{gather*}
Of course, ASM having the same f\/irst (resp. last) $k$ rows
correspond to staircases having the same  $k$ right (resp. left)
columns. Moreover, given a staircase in the letters $1,\ldots,n$,
one can complete it in a canonical manner into a staircase with
columns of lengths $n,\ldots,1$. The partition function of a set
of staircases having two f\/ixed columns $u,v$ factorizes into
$F(n,u) F(u,v) F(v,[ \,])$. The middle part can be interpreted as
the partition function (Theorem \ref{th:Esc2Col}) of a set of ASM
having a f\/ixed top and bottom, up to a factor due to the f\/ixed
top and bottom parts.

\medskip
Instead of staircases, one can use \emph{ribbon tableaux} (a
\emph{ribbon} is a skew diagram which does not contain a
$2\times2$ sub-diagram, a \emph{ribbon tableau} is an increasing
sequence of diagrams of partitions which dif\/fer by a ribbon).
Indeed, given $n$ and a column $u$ of length $n-k$ such that
$u_1\leq n$, let $\widetilde u = \{1,\ldots,n\} \setminus u$
(sorted increasingly). Def\/ine
\begin{gather*}
\p(u,n) := [\widetilde u_1-1, \ldots, \widetilde u_k-k].
\end{gather*}

It is immediate to translate the condition that a sequence of
columns
       is a staircase in terms of diagrams of partitions.

\begin{lemma}
Given $n$, a column $u$, $k=n-\ell(u)$. Then the map
\begin{gather*}
v\to \mu= \p(v,n)
\end{gather*}
is a bijection between the set of columns $v$, $v_1\leq n$,
 such that $vu$ is  a staircase and the set of ribbons
$ \bigl(\p(u,n)+1^k \bigr)/\mu$, $\mu\in \N^{k-1}$.
\end{lemma}

For example, for $u=[5,3,2]$, $n=6$, then
\begin{gather*}
\p(u,6)=[1\moins 1, 4\moins 2, 6\moins 3]= [0,2,3],
\end{gather*}
hence
 $\p(u,6)+1^3= [1,3,4]$. There are eleven ribbons $[1,3,4]/\mu$, $\mu\in\N^2$,
in bijection with the staircases~$[v,u]$:
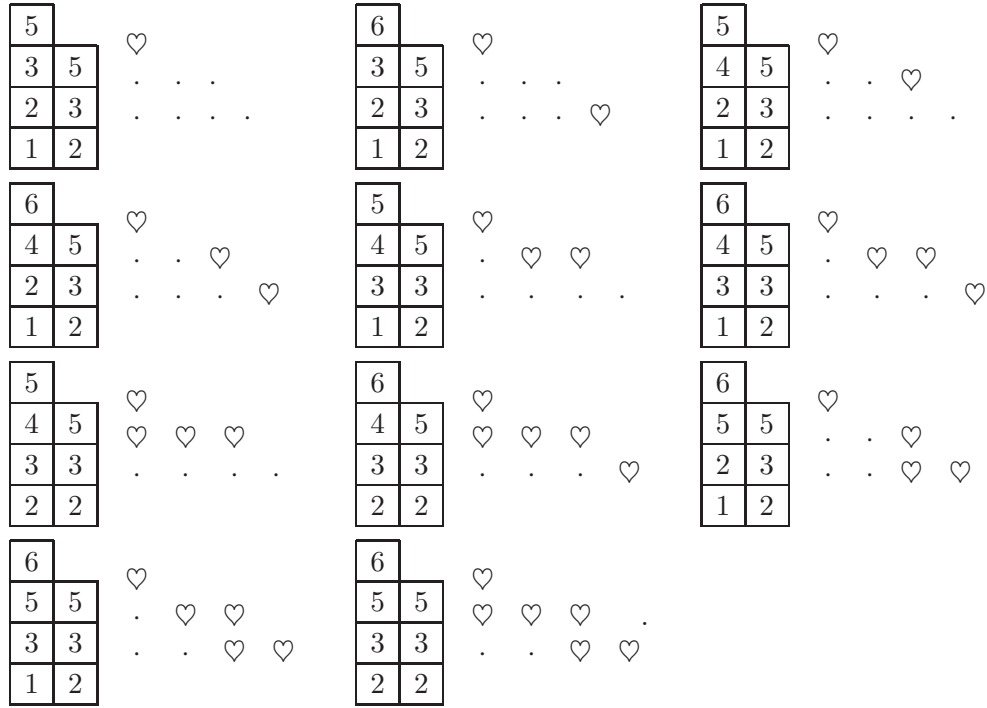
\begin{gather*}
\begin{array}{lll}
\young{
 5 \cr
 3 & 5 \cr
 2 & 3 \cr
 1 & 2 \cr
}\quad
\begin{matrix}
\heartsuit\cr \cdot &\cdot &\cdot \cr \cdot &\cdot &\cdot &\cdot
\cr
\end{matrix} \, \quad
&\young{
 6 \cr
 3 & 5 \cr
 2 & 3 \cr
 1 & 2 \cr
}\quad
\begin{matrix}
\heartsuit\cr \cdot &\cdot &\cdot \cr \cdot &\cdot &\cdot
&\heartsuit\cr
\end{matrix} \, \quad
&\young{
 5 \cr
 4 & 5 \cr
 2 & 3 \cr
 1 & 2 \cr
} \quad
\begin{matrix}
\heartsuit\cr \cdot &\cdot &\heartsuit\cr \cdot &\cdot &\cdot
&\cdot \cr
\end{matrix}
\\[6ex]
\young{
 6 \cr
 4 & 5 \cr
 2 & 3 \cr
 1 & 2 \cr
}\quad
\begin{matrix}
\heartsuit\cr \cdot &\cdot &\heartsuit\cr \cdot &\cdot &\cdot
&\heartsuit\cr
\end{matrix} \, \quad
&\young{
 5 \cr
 4 & 5 \cr
 3 & 3 \cr
 1 & 2 \cr
}\quad
\begin{matrix}
\heartsuit\cr \cdot &\heartsuit &\heartsuit\cr \cdot &\cdot &\cdot
&\cdot \cr
\end{matrix}\, \quad
&\young{
 6 \cr
 4 & 5 \cr
 3 & 3 \cr
 1 & 2 \cr
}\quad
\begin{matrix}
\heartsuit\cr \cdot &\heartsuit &\heartsuit\cr \cdot &\cdot &\cdot
&\heartsuit\cr
\end{matrix}
\\[6ex]
 \young{
 5 \cr
 4 & 5 \cr
 3 & 3 \cr
 2 & 2 \cr
}\quad
\begin{matrix}
\heartsuit\cr \heartsuit &\heartsuit &\heartsuit\cr \cdot &\cdot
&\cdot &\cdot \cr
\end{matrix}\, \quad
&\young{
 6 \cr
 4 & 5 \cr
 3 & 3 \cr
 2 & 2 \cr
}\quad
\begin{matrix}
\heartsuit\cr \heartsuit &\heartsuit &\heartsuit\cr \cdot &\cdot
&\cdot &\heartsuit\cr
\end{matrix}\, \quad
&\young{
 6 \cr
 5 & 5 \cr
 2 & 3 \cr
 1 & 2 \cr
}\quad
\begin{matrix}
\heartsuit\cr \cdot &\cdot &\heartsuit\cr \cdot &\cdot &\heartsuit
&\heartsuit\cr
\end{matrix}\, \quad
\\[6ex]
\young{
 6 \cr
 5 & 5 \cr
 3 & 3 \cr
 1 & 2 \cr
}\quad
\begin{matrix}
\heartsuit\cr \cdot &\heartsuit &\heartsuit\cr \cdot &\cdot
&\heartsuit &\heartsuit\cr
\end{matrix}\, \quad
&\young{
 6 \cr
 5 & 5 \cr
 3 & 3 \cr
 2 & 2 \cr
}\quad
\begin{matrix}
\heartsuit\cr \heartsuit &\heartsuit &\heartsuit\cr \cdot &\cdot
&\heartsuit &\heartsuit\cr
\end{matrix}.
\end{array}
\end{gather*}

\vspace{13pt} We shall now give a weight to a ribbon $\zeta/\mu$.
Number boxes of the diagram of $\zeta$ uniformly in each diagonal,
by $1,2,3,\ldots$, starting from the top leftmost box. For each
box  $\square$ of $\zeta/\mu$, denote $c(\square)$ this number
(this is a shifted \emph{content} \cite[p.~11]{Mac}). A box of a
ribbon is \emph{terminal} if it is the rightmost in its row.

Given two partitions $\zeta,\mu$ such that $\zeta/\mu$ is a
ribbon, we shall weigh the ribbon by giving a~weight to each of
its boxes as follows:
\begin{itemize}\itemsep=0pt
\item  a box  $\square$  which is not terminal is weighted
$x-y_{c(\square)}$; \item a box which is terminal has weight
$y_{c(\square)}$ if it is above another box, or weight $x$ if not;
\item $\theta(\zeta/\mu)$ is the product of these elementary
weights (this is a polynomial in $x,y_1,y_2,\ldots$).
\end{itemize}

For example, for $\zeta=[3,3,3,5,5]$, $\mu=[1,2,2,3,5]$, one has
the following weights:
\begin{gather*}
\begin{matrix}
 \cdot & x\moins y_2  & y_3 \\
\cdot & \cdot & y_4 \\
\cdot & \cdot & x \\
\cdot & \cdot & \cdot & x\moins y_7  & x \\
\cdot & \cdot & \cdot &\cdot & \cdot
\end{matrix}
\end{gather*}

The link between the weight of ribbons and the weight of
staircases will appear in the proof of the next theorem.

\section{Right truncated staircases}

Given a column $u=[u_1,\ldots, u_r]$, write the numbers
$u_1,\ldots,1$ inside a ribbon, passing to a new level for each
value belonging to $u$, $u_1$ being at level $0$. The sequence of
levels  of $1,2,\ldots, u_1$ is denoted $\langle u \rangle$.
\begin{gather*}
u=[\mathbf{5}, \mathbf{3}, \mathbf{2}]\quad  \Rightarrow \young{
\mathbf{5} \cr 4 & \mathbf{3}\cr \blk &\mathbf{2} \cr \blk &1 \cr}
\quad  \Rightarrow \young{ 0 \cr 1 & 1\cr \blk &2 \cr \blk &3 \cr}
\quad  \Rightarrow\  \langle u \rangle = [3,2,1,1,0].
\end{gather*}

\begin{theorem}
Let $n$ be a positive integer, let $u$ be a column of length
$n-r$, such that $u_1\leq n$. Then $F(n,u) = x^{\rho_r}
y^{-\langle  \widetilde u \rangle} \X_{\widetilde u, u^\omega}(
\x, \y)$, with $u^\omega=[u_{n-r},\ldots, u_1]$, $\widetilde{u}$
the (increasing) complement of~$u$ in $[1,\ldots,n]$, and
$\rho_r=[r\moins 1,\ldots,1,0]$.
\end{theorem}

\begin{proof}
We shall decompose the set $\cE(n,u)$ according to the
ante-penultimate column $v$, assuming the theorem true for columns
of length $n-r+1$. This translates into the equality
\begin{gather*}
   F(n,u;\x^r) = \sum_v  F(n,v;\x^{r-1})   F(v,u; x),
\end{gather*}
where the sum is over all columns $v$ of length $n-r+1$, writing
$x$ for $x_r$.

In terms of Schubert polynomials, one therefore is reduced to show
that
\begin{gather*}
  (x_1\cdots x_{r-1}) y^{-\langle \widetilde u\rangle }
   \X_{\widetilde u, u^\omega}( \x^{r}, \y)    =
 \sum_v  y^{-\langle \widetilde v\rangle }
   \X_{\widetilde v, u^\omega}( \x^{r-1}, \y)
    F(v,u;x).
\end{gather*}

Proposition \ref{th:ProdSchub} gives the expansion of $(x_1\cdots
x_r) \X_{\widetilde u, u^\omega}( \x^{r}, \y)$, the
coef\/f\/icients being weights of ribbons. We have therefore to
compare these weights, multiplied by a monomial in $y$, and
divided by $x=x_r$, to the appropriate $F(v,u;x)$ to be able to
conclude.

In more details, given $v,u$, then $\zeta=[\widetilde u_1,
\widetilde u_2 -1,\ldots,
 \widetilde u_{n-r} -n+r+1]$,
$\mu= [\widetilde v_1-1, \ldots, \widetilde v_{n-r-1} -n+r+1]$,
and the ribbon is $\zeta/\mu$. One notices that $ y^{\langle
v\rangle }/y^{\langle u\rangle } = y^{\langle \widetilde u\rangle
}/y^{\langle \widetilde v\rangle }$ is a monomial without
multiplicity.  The $i$'s such that $y_i$ appear in it are the
contents of the boxes which are not terminal boxes at the bottom
of their column. But these values are exactly those which
correspond to a~factor $x/y_i$ or $(x/y_i -1)$ in $F(v,u;x)$.
Moreover, it is clear that the contents of boxes which are not
terminal boxes are the numbers which appear in $u$ and $v$ at the
same level. Therefore the expansion of the Schubert polynomial
$\X_{\widetilde u, u^\omega}( \x^{r}, \y)$, multiplied by
$x_1\cdots x_{r-1}  x$, and the expansion of $F(n,u)$ according to
the last variable $x=x_r$ coincide, once normalized, and the
theorem is proved.
\end{proof}

For example, for $v=[16,15,12,10,9,5,3,2,1]$,
$u=[15,14,10,9,5,4,3,1]$, and $n=17$, one has $\widetilde v=
[4,6,7,8,11,13,14,17]$, $\mu=[3,4,4,4,6,7,7,9]$,
$\widetilde u=[2,6,7,8,11,12,13,17,17]$,\\
$\zeta=[2,5,5,5,7,7,7,9,9]$. Moreover,
$\langle v\rangle =[8, 7, 6, 6, 5, 5, 5, 5, 4, 3, 3, 2, 2, 2, 1]$,\\
$\langle u\rangle =[7, 7, 6, 5, 4, 4, 4, 4, 3, 2, 2, 2, 2, 1]$ and
$y^{\langle v\rangle  -\langle u\rangle}=
y_{1}y_{4}y_{5}y_{6}y_{7}y_{8}y_{9}y_{{10}}y_{{11}}y_{{14}}y_{15}$.
\begin{gather*}
\young{ 16 \cr 15 &15 \cr 12 &14 \cr 10 &10 \cr
 9 & 9 \cr
 5 & 5 \cr
 3 & 4 \cr
 2 & 3 \cr
 1 & 1 \cr
} \qquad
\begin{matrix}
\fbox{$1$} & x\cr \cdot &\cdot &\cdot &\fbox{$5$} & 6 &\cr \cdot
&\cdot &\cdot &\cdot &7 &\cr \cdot &\cdot &\cdot &\cdot &8 &\cr
\cdot &\cdot &\cdot &\cdot &\fbox{$9$} &\fbox{$10$} &11 &\cr \cdot
&\cdot &\cdot &\cdot &\cdot &\cdot &x &\cr \cdot &\cdot &\cdot
&\cdot &\cdot &\cdot &\cdot &\cr \cdot &\cdot &\cdot &\cdot &\cdot
&\cdot &\cdot &\fbox{$15$} & x &\cr \cdot &\cdot &\cdot &\cdot
&\cdot &\cdot &\cdot &\cdot &\cdot &\cr
\end{matrix}
\end{gather*}

The weights are represented, writing $i$ for $x y_i^{-1}$ and
\fbox{$ i$} for $x y_i^{-1} -1$.

As a special case of the theorem, one can f\/ilter the complete
staircases according to their column of length $1$.

\begin{corollary}
Let  $n$, $b$ be two positive integers, $b\leq n$. Then
\begin{gather*}
   F(n, [b])=  x^\rho  y^{-\rho} y^{-[0^{b-1}, 1^{n-b}]}
 \Y_{[0^{b-1}, 1^{n-b}]}(\x,\y),
\end{gather*}
with $\rho=[n-2,\ldots,0]$, $\x=[x_1, \ldots, x_{n-1}]$.
\end{corollary}

The Schubert polynomials appearing in the corollary specialize,
for $\y=\mathbf{0}$, to the elementary symmetric functions in
$x_1,\ldots, x_{n-1}$.

For example, for $n=3$, 2 staircases contribute to $F(3,[1])=
x_1/y_1 y^{-[1,1]} \Y_{11}(\x,\y)$, 3 staircases contribute to
$F(3,[2])= x_1/y_1 y^{-[0,1]} \Y_{01}(\x,\y)$, and the last two
contribute to $F(3,[3])= x_1/y_1 y^{-[0,0]} \Y_{00}(\x,\y)$.

\medskip
In our opinion, the most fundamental property shown by the above
theorem is the symmetry, in the variables $x_i$'s, of the function
$F(n,u) x^{-\rho}$.

\section{General staircases}

Let $n$, $k$, $r$ be three positive integers, and $v$ be a column
of length $n-k-r \geq 0$. Put $\x=[x_1,\ldots, x_k]$,
$\z=[x_{k+1},\ldots, x_{k+r}]$. Then
\begin{gather*}
  F(n,v;\x, \z )
=  x^{\rho_k} (x_1\cdots x_k)^r z^{\rho_r}
 y^{-\langle\widetilde u\langle}   \X_{\widetilde u, u^\omega}( \x, \y)
 =  \sum_u F(n,u;\x)   F(u,v;\z),
\end{gather*}
sum over all columns $u$ of length $n-k$.

Since Schubert polynomials in $\x$ are linearly independent (and
so are their products by a~f\/i\-xed~$x^\rho$), then
 $F(u,v;\z)$ can be characterized as the coef\/f\/icient
of $\X_{\widetilde u, u^\omega}( \x, \y) $ in
\begin{gather*}
  x^{-\rho_k} F(n,v;\x,\z )  y^{\langle\widetilde v\rangle}= (x_1\cdots x_k)^r  z^{\rho_r}  \X_{\widetilde v, v^\omega}
( \x,\z ;  \y).
 \end{gather*}

Therefore it is the specialization $x_1=y_1, x_2=y_2,\ldots$ of
the image of this polynomial under $\partial_\sigma$, with
$\sigma= [\widetilde u, u^\omega]$, thanks to Lemma
\ref{th:Newton}. Let $\beta=[\widetilde v_1 -1,\ldots, \widetilde
v_{k+r} -1]$,
 $\alpha= [\beta_1+r,\ldots, \beta_{k+r}-k+1]$.
Writing the Schubert polynomial $\X_{\widetilde v, v^\omega}$
 as a determinant
\begin{gather*}
  S_{\beta_1-0,\ldots, \beta_{r+k}-r-k+1}
(\x+\z -\y^{\beta_1},\ldots, \x+\z-\y^{\beta_{r+k}}),
\end{gather*}
then its product by $x_1^r\cdots x_k^r z_1^r\cdots z_r^r$ is equal
to
\begin{gather*}
  S_{\alpha} (\x+\z -\y^{\beta_1},\ldots, \x+\z-\y^{\beta_{r+k}}).
\end{gather*}
Thanks to (\ref{SchurDrap}), the value of this determinant is not
changed by replacing, in rows $1,2,\ldots, n\plus k$,  $\x=\x^k$
by $\mathbf{0},\ldots, \mathbf{0}, \x^1,\ldots, \x^k$
respectively.

Now, one can easily compute the image of this last determinant
under $\partial_\sigma$. The transformation boils down to
decreasing, in each row the indices of complete functions by some
quantity, and increasing at the same time the exponent of $\x^j$
by the same amount.

The resulting determinant  is
\begin{gather*}
   S_{\alpha/[0^r,\gamma]} (\z -\y^{\beta_1},\ldots, \z-\y^{\beta_{r+k}} //
 \mathbf{0},\ldots, \mathbf{0}, \x^{1+\gamma_1},\ldots, \x^{k+\gamma_k}   ),
\end{gather*}
with $\gamma=[\tilde u_1-1,\ldots, \tilde u_k-k]$. The last
transformation is to replace $\x$ by $\y$, and this produces
\begin{gather*}
   S_{\alpha/[0^r,\gamma]}  (\z -\y^{\beta_1},\ldots, \z-\y^{\beta_{r+k}} //
 \mathbf{0},\ldots, \mathbf{0}, \y^{1+\gamma_1},\ldots, \y^{k+\gamma_k}   ).
\end{gather*}

In conclusion, one has the following theorem.

\begin{theorem}    \label{th:Esc2Col}
Let $n$, $k$, $r$ be three positive integers,  $v$ be a column of
length $n-k-r$, $u$ be a column of length $n-k$, $\z=[z_1,\ldots,
z_r]$. Let $\beta=[\widetilde v_1 -1,\ldots, \widetilde v_{k+r}
-1]$, $\alpha= [\beta_1+r,\ldots, \beta_{k+r}-k+1]$,
$\gamma=[\tilde u_1-1,\ldots, \tilde u_k-k]$. Then
\begin{gather}
     F(u,v;\z) = z^{\rho_r}   y^{\langle v\rangle -\langle u\rangle}
  S_{\alpha/[0^r,\gamma]}  (\z -\y^{\beta_1},\ldots, \z-\y^{\beta_{r+k}} //
 \mathbf{0},\ldots, \mathbf{0}, \y^{1+\gamma_1},\ldots, \y^{k+\gamma_k}).
\label{MultiSchur}
\end{gather}
\end{theorem}

For example, for $n=6$, $u=[6,5,3,1]$, $v=[5,2]$, one has
$\widetilde v =[1346]$, $\widetilde u=[24]$, $\beta=[0,2,3,5]$,
$\alpha=[2,3,3,4]$, $\gamma=[1,2]$,
\begin{gather*}
 F(6,[5,2]) = x^{3210}/y^{33211} \X_{1346 25}  = x^{3210}/y^{33211} \Y_{011200},
\\
 F((6,[6,5,3,1]) = x^{10}/ y^{211} \X_{24 1356}  = x^{10}/ y^{211} \Y_{12 0000},
\\
\Y_{011200} = S_{0112}(\x^4-\y^0, \x^4-\y^2, \x^4-\y^3,\x^4-\y^5)
,
\\
 x^{2222} \Y_{011200} = S_{2334}(\x^4-\y^0, \x^4-\y^2, \x^4-\y^3,\x^4-\y^5).
\end{gather*}
Putting $[x_3,x_4]=\z$, the last determinant is written
\begin{gather*}
 S_{2334}(\z-\y^0, \z-\y^2, \z-\y^3,\z-\y^5  //  \mathbf{0},\mathbf{0},
  \x^1, \x^2).
\end{gather*}
Its image under $\partial_\sigma = \partial_{24 1356}= \partial_1
\big(\partial_3(\partial_2(  ))\big)$ is
\begin{gather*}
 S_{2334/0012}(\z-\y^0, \z-\y^2, \z-\y^3,\z-\y^5  //
         \mathbf{0},\mathbf{0}, \x^{1+1}, \x^{2+2}).
\end{gather*}
In f\/inal,
\begin{gather*}
F([6,5,3,1], [5,2])
\\ \qquad
{}= x^{10} y^{211}/y^{33211}\cdot  \begin{vmatrix}
S_2(\z )  & S_4(\z -\y^2 ) & S_5(\z-\y^3 ) &S_7(\z-\y^5 ) \\
S_1(\z )  & S_3(\z-\y^2 ) & S_4(\z-\y^3 ) &S_6(\z-\y^5 )    \\
0         & S_1(\z-\y^2+\y^2 ) & S_2(\z-\y^3+\y^2) &S_4(\z-\y^5+\y^2 )    \\
0         & 0                  & S_0(\z-\y^3+\y^4 )
&S_2(\z-\y^5+\y^4)
\end{vmatrix}.
\end{gather*}

\medskip
Notice that the expression of $F(u,v;\z) z^{-\rho_r}$ given in
(\ref{MultiSchur}) is symmetrical in $z_1,\ldots,z_r$. It would be
interesting to prove directly that $F(u,v;z_1,z_2) z_1^{-1}$, for
$\ell(u)=\ell(v)+2$, is a symmetric function in $z_1$, $z_2$.

\section{Left truncated staircases}

Let us now treat the staircases $[u=u_1,u_2,\ldots, u_{r+1}=[
\,]]$, with a f\/ixed left column $u$. Let $n$ be such that
$u_1\leq n$.

{\sloppy
>From the preceding section, one knows that
$F(u,[ \,];z_1,\ldots,z_r) z^{-\rho_r}$ is the  coef\/f\/icient of
$F(n,u) x^{-\rho_{n\!-\! r}}$ in the expansion of
\begin{gather*}   x_1^r\cdots x_{n-r}^r
          = \X_{[r+1,\ldots,n,1,\ldots,r]}(\x, \mathbf{0}).
\end{gather*}

}

This expansion is a special case of Cauchy formula \cite[Theorem
10.2.6]{CBMS} for three alphabets, and any permutation $\sigma$:
\begin{gather*}
 \X_\sigma(\x, \mathbf{w}) = \sum_{\sigma',\sigma''}
 \X_{\sigma'}(\y, \mathbf{w})  \X_{\sigma''}(\x, \y),
\end{gather*}
sum over all reduced products $\sigma= \sigma'  \sigma''$, in the
case where $\sigma$ is a Gra{\ss}mannian permutation, and
$\mathbf{w}= \mathbf{0}:=[0,0,\ldots]$.

Using that $\X_\tau(\y,\mathbf{0})= \X_{\tau^{-1}}(\mathbf{0},
  \overline{\y})$, for any permutation $\tau$,
with $\overline{\y} = [-y_1,-y_2,\ldots]$, going back to the $x$
variables instead of $z$, and taking into account the value of
$u_r$ (which introduces a shift of indices, compared to the case
$u_r=1$), one has the following theorem.

\begin{theorem}
Let $u$ be a column of length $r$, with $u_r=1$. Then
\begin{gather*}
F(u,[ \,]) = x^{\rho_r}  y^{-\langle u\rangle}
  \X_{u^\omega, \widetilde{u}}(\mathbf{0}, \overline{\y}),
\end{gather*}
with $\overline{\y}=[-y_1,-y_2, \ldots]$, $u^\omega=[u_r,\ldots,
u_1]$, $\widetilde{u}$ the (increasing) complement of~$u$ in
$[1,\ldots,n]$, and~$\langle u\rangle$  defined Section~$4$.

If $u_r=k>1$, then $F(u,[ \,])$ is obtained from $F([u_r\moins
k\plus1, \ldots, u_1\moins k\plus 1],  [ \,])$ by increasing the
indices of all $y_i$ by $k-1$.
\end{theorem}

For example, if $u=[5,3,1]$, $n=6$, then $\langle
u\rangle=[2,2,1,1]$
\begin{gather*}
 F([5,3,1],[])= x^{210}  y^{-2,-2,-1,-1}
\X_{135246}(\mathbf{0}, \overline{\y}).
\end{gather*}
If  $u=[6,4,2]$,  then
\begin{gather*}
 F([6,4,2],[])= x^{210}  y^{0,-2,-2,-1,-1}
 \X_{135246}(\mathbf{0}, [-y_2,-y_3,\ldots]).
\end{gather*}

\section{Appendix}

We need some combinatorial properties of Schubert polynomials
(Lemma \ref{th:VerticalStrip}, Proposition \ref{th:ProdSchub}),
completing those which can be found in \cite{FK, CBMS,MacSchub}.

Schubert polynomials $\X_\sigma(\x,\y)$ constitute a linear basis
of the ring of polynomials in $x_1,x_2,\ldots$, with
coef\/f\/icients in $y_1,y_2,\ldots$, indexed by permutations in
$\mfS_\infty$ (permutations f\/ixing all, but a~f\/inite number of
integers). Their family is stable under divided dif\/ferences, and
they all vanish when $\x$ is specialized to $\y$ (i.e.\ $x_1\to
y_1$, $x_2\to y_2,\ldots$), except the polynomial
$\X_{1}(\x,\y)=1$. In fact, these two properties, added to some
normalizations, characterize Schubert polynomials uniquely.

This allows, for example, to expand any polynomial in $\x$ in the
basis of Schubert polynomials \cite[Theorem~9.6.1]{CBMS}:
\begin{lemma}   \label{th:Newton}
Given a polynomial $f(x_1,x_2,\ldots)$, then
\begin{gather*}
  f(\x)= \sum_{\sigma\in \mfS_\infty} \partial_{\sigma}
        \bigl(f(\x)\bigr)\bigr|_{\x=\y}  \X_{\sigma}(\x,\y).
\end{gather*}
\end{lemma}

\begin{proof}
 Take the image of both members under some $\partial_\sigma$.
The right hand-side is still a sum of Schubert polynomials, with
modif\/ied indices.  Specialize now $\x$ to $\y$. There is only
one term which survives, the Schubert polynomial $\X_1$, which
comes from $\X_{\sigma}(\x,\y)$. Schubert polynomials being a
linear basis, these equations for all $\sigma$ determine the
function $f(x)$.
\end{proof}

Schubert polynomials can be interpreted in terms of tableaux.
 This gives a simple way of obtaining the branching
rule of these polynomials according to the last variable $x_r$.

In the special case of a Gra{\ss}mannian polynomial, with descent
in $r$, using increasing partitions, the branching is
\begin{gather}
\label{HopfSchub}
  \Y_{\lambda}(\x^r,\y)  = \sum_\mu
  \psi^h\biggl(\frac{\lambda}{\mu}; r\biggr)  \Y_{\mu}(\x^{r-1},\y)
\end{gather}
sum over all $\mu$ such that $\lambda/\mu$ be an horizontal strip,
the weight of the strip being
\begin{gather*}
 \psi^h(\lambda/\mu; r) :=\prod_{\square \in \lambda/\mu}
  (x_r - y_{r+c(\square)}),
\end{gather*}
product over all boxes of the strip, $c(\square)$ being the usual
content, i.e.\ the distance of $\square$ to the main diagonal of
the diagram of $\lambda$.

Given two partitions $\zeta$, $\mu$ such that $\zeta/\mu$ be a
vertical strip, def\/ine similarly the \emph{weight} of the strip
to be
\begin{gather*}
 \psi^v\biggl(\frac{\zeta}{\mu};r\biggr) = \prod_{\square \in \zeta/\mu}  y_{r+c(\square)}.
\end{gather*}

This weight appear in the product of a Gra{\ss}mannian polynomial
by a monomial, as states the following lemma.

\begin{lemma}    \label{th:VerticalStrip}
Let $\nu\in\N^r$ be a partition, $\zeta=\nu+1^r$. Then
\begin{gather*}
 x_1\cdots x_r  \Y_\nu(\x, \y) = \sum_\mu \psi^v\biggl(\frac{\zeta}{\mu};r\biggr)  \Y_\mu(\x,\y),
\end{gather*}
sum over all $\mu$ such that $\zeta/\mu$ be a vertical strip.
\end{lemma}

\begin{proof}
  Let $\lambda= [\nu_r+r-1,\ldots, \nu_1+0]$.
 The polynomial $\Y_\nu(\x,\y)$ is by def\/inition the image of
\begin{gather*}
 \Y_{\lambda}(\x,\y) = \prod_{(i,j)\in\lambda} (x_i-y_j)
\end{gather*}
under $\partial_{\omega_r}$.  Since multiplication by $x_1\cdots
x_r$ commutes with $\partial_{\omega_r}$, then $x_1\cdots x_r
\Y_\nu(\x, \y)$  is the image of
\begin{gather*}
 x_1\cdots x_r\Y_{\lambda}(\x,\y)  =
\bigl( (x_1-y_{\nu_r+r}) + y_{\nu_r+r}\bigr)  \cdots
 \bigl( (x_r-y_{\nu_1+1}) + y_{\nu_1+1}\bigr)  \Y_\nu(\x, \y).
\end{gather*}
This polynomial expands into a sum of Schubert polynomials,
multiplied by some $y_j$. Taking the image of this sum under
$\partial_{\omega_r}$ involves only transforming the indices.
Explicitly,
\begin{gather*}
 x_1\cdots x_r  \Y_\nu(\x, \y) =
  \sum_{\epsilon\in [0,1]^n}
   y_{\nu_1+1}^{\epsilon_1}\cdots y_{\nu_r+r}^{\epsilon_r}
  \Y_{\zeta- \epsilon},
\end{gather*}
sum over all $\epsilon$ such that $\zeta- \epsilon$ be a
partition, say $\mu$. In other words, $\zeta/\mu$ is  a vertical
strip and the lemma is proved.
\end{proof}

For example, for $\nu= [2,2,3]$,
\begin{gather*}
 x_1x_2x_3 \Y_{223} = \sum_{\epsilon}
 y_3^{\epsilon_1} y_4^{\epsilon_2} y_6^{\epsilon_3}  \Y_{334-\epsilon}
 \\ \phantom{x_1x_2x_3 \Y_{223}}
{}= \Y_{334} +y_3 \Y_{234} +y_6 \Y_{333} + y_3y_4 \Y_{224} +y_3y_6
\Y_{233} +y_3y_4y_6 \Y_{223},
\end{gather*}
discarding the terms $y_4 \Y_{324}$, $y_4y_6 \Y_{323}$ as being
non-conform.

\begin{proposition}   \label{th:ProdSchub}
Let $\nu\in\N^r$ be a partition, $\zeta=\nu+1^r$. Then
\begin{gather*}
x_1\cdots x_r \Y_\nu(\x^r,\y) \bigr|_{x_r=x} = \sum_\mu
\theta\biggl(\frac{\zeta}{\mu}\biggr)   \Y_{\mu}(\x^{r-1},\y),
\end{gather*}
sum over all partitions $\mu\in\N^{r-1}$ such that $\zeta/\mu$ be
a ribbon, $\theta(\zeta/\mu)$ being defined in Section~$3$.
\end{proposition}

\begin{proof}
The polynomial $x_1\cdots x_r \Y_\nu(\x^r,\y)$ is a sum of
Schubert polynomials obtained by suppressing vertical strips to
$\zeta$, according to Lemma \ref{th:VerticalStrip}. Expanding then
according to $x_r$ involves removing horizontal strips. The
resulting partitions $\mu$ dif\/fer from $\zeta$ by a ribbon.
Given such a~$\mu$, the coef\/f\/icient of $\Y_\mu$ in the f\/inal
polynomial will the sum
\begin{gather}
\label{DecompoRuban}
 \sum_\eta \psi^v\biggl(\frac{\zeta}{\eta};r\biggr)  \psi^h\biggl(\frac{\eta}{\mu}\biggr)
\end{gather}
sum over all $\eta$ such that $\zeta/\eta$ be a vertical strip and
$\eta/\mu$   be an horizontal strip.

The non-terminal boxes of $\zeta/\mu$ are common to all horizontal
strips, the terminal boxes which are above another box of the
ribbon are common to all vertical strips. Hence, up to a common
factor, the sum (\ref{DecompoRuban}) reduces to
\begin{gather*}
  \sum_{A'\cup A''} \prod_{\square\in A'} y_{c(\square)}
  \prod_{\square\in A''}(x- y_{c(\square)}),
\end{gather*}
sum over all decompositions of the remaining set $A$ of boxes of
$\zeta/\mu$ into two disjoint subsets. This sum is clearly equal
to $x^k$, $k$ being the cardinality of $A$. In total, the factor
of $Y_\mu(\x^{r-1},\y)$ in the expansion  of $x_1\cdots x_r
\Y_\nu(\x^r,\y)$ is precisely $\theta(\zeta/\mu)$.
\end{proof}

\subsection*{Acknowledgements}
The author benef\/its from the ANR project BLAN06-2\_134516.

\pdfbookmark[1]{References}{ref}
\LastPageEnding

\end{document}